\documentclass[11pt]{amsart}
\usepackage{tabularx,booktabs}
\usepackage{caption}
\usepackage{amsmath}
\usepackage{amsfonts}

\usepackage{amscd}
\usepackage{amsthm}
\usepackage{amssymb} \usepackage{latexsym}
\usepackage{eufrak}
\usepackage{euscript}
\usepackage{epsfig}
\usepackage{graphics}
\usepackage{array}
\usepackage{enumerate}
\usepackage{dsfont}
\usepackage{color}
\usepackage{wasysym}
\usepackage{hyperref}
\usepackage{pdfsync}

\newcommand{\bel}[1]{\begin{equation}\label{#1}}

\newcommand{\be}{\begin{equation}}

\newcommand{\ba}{\begin{eqnarray}}
\newcommand{\ea}{\end{eqnarray}}

\newcommand{\qe}{\end{equation}}

\newcommand{\Hmm}[1]{\leavevmode{\marginpar{\tiny%
$\hbox to 0mm{\hspace*{-0.5mm}$\leftarrow$\hss}%
\vcenter{\vrule depth 0.1mm height 0.1mm width \the\marginparwidth}%
\hbox to
0mm{\hss$\rightarrow$\hspace*{-0.5mm}}$\\\relax\raggedright #1}}}

\newtheorem{theorem}{Theorem}[section]

\newtheorem{lemma}[theorem]{Lemma}
\newtheorem{corollary}[theorem]{Corollary}

\newtheorem{remark}[theorem]{Remark}

\newtheorem{prop}[theorem]{Proposition}

\makeatletter

\@addtoreset{equation}{section}

\makeatother
\begin{document}

\title[A parabolic flow approach to generalized Gauduchon metrics]{A parabolic flow approach to generalized Gauduchon metrics}

\author{Han Luo}
\address{Han Luo: School of Mathematics and Statistics, Hainan University, Haikou
570228, Hainan, People’s Republic of China }
\email{\href{mailto:luohan@hainanu.edu.cn}{luohan@hainanu.edu.cn}}

\begin{abstract}
We give a parabolic proof of the Fu–Wang–Wu theorem for semilinear elliptic equations associated with generalized Gauduchon metrics. Fixing a Hermitian metric $\omega$, we introduce a flow $\partial_t v-\Delta v=\psi(\vert\nabla v\vert^2)+\langle B,dv\rangle-c(t)-f(x)$ for the conformal factor, whose stationary points solve the semilinear elliptic equation of Fu--Wang--Wu. Under the growth condition $\liminf_{t\to+\infty}\psi(t)/t^\mu>0$ for some $\mu>\tfrac{1}{2}$, we prove that the flow exists uniquely for all time and converges; the limit solves the elliptic equation and is unique up to addition of a constant. Our convergence argument does not rely on a gradient-flow structure or a monotone energy functional. In particular, this recovers, by a flow method rather than the continuity method, the existence and uniqueness theorem of Fu--Wang--Wu on $k$-th generalized Gauduchon metrics, providing a new parabolic approach to their elliptic theory.
\end{abstract}

\maketitle

\section{Introduction}

Let $(M,\omega)$ be a compact Hermitian manifold of complex dimension $n$, without boundary. Special Hermitian metrics play an important role in non-K\"ahler geometry. For $n\geq 2$, Gauduchon's theorem \cite{Gau77,Gau84} states that every conformal class contains a metric satisfying $\partial\bar\partial\omega^{n-1}=0$, unique up to multiplication by a positive constant.

Fu--Wang--Wu \cite{FWW13} generalized this condition as follows. A metric $\omega$ is called a $k$-th generalized Gauduchon metric, or $k$-Gauduchon metric, if
\begin{equation}\label{eq:intro-k-gauduchon}
    \partial\bar\partial\omega^k\wedge\omega^{n-k-1}=0,
    \qquad 1\leq k\leq n-1.
\end{equation}
This recovers the classical Gauduchon condition when $k=n-1$. Pluriclosed (SKT) metrics (see \cite{ST10}) satisfy the condition for $k=1$, while astheno-K\"ahler metrics in the sense of Jost--Yau \cite{JY93,JY94} satisfy it for $k=n-2$ when $n\geq 3$. Examples and related structures have been studied in \cite{FU13,LU17,LUV17,FGV19}.

Writing $\widetilde\omega=e^{v/k}\omega$, the $k$-Gauduchon condition becomes
\begin{equation*}
    \Delta v+|\nabla v|^2+\langle B,dv\rangle=-\varphi,
    \qquad
    \varphi=
    \frac{n(\sqrt{-1}/2)\partial\bar\partial(\omega^k)
    \wedge\omega^{n-k-1}}{\omega^n},
\end{equation*}
where $B$ is a smooth real $1$-form determined by $\omega$ and $k$, and $\Delta$ and $\nabla$ are associated with $\omega$. This equation need not be solvable for a given background metric. Fu--Wang--Wu therefore introduce an additional unknown constant and, more generally, establish by the method of continuity an existence and uniqueness theorem for
\begin{equation}\label{eq:intro-elliptic}
    \Delta v+\psi(|\nabla v|^2)+\langle B,dv\rangle=c+f(x),
\end{equation}
where $\psi\in C^\infty(\mathbb R)$ satisfies a growth condition, $B$ is a smooth real $1$-form, and $f\in C^\infty(M)$ has zero mean with respect to $\omega^n$.

Parabolic methods provide another approach to elliptic problems in complex geometry, as illustrated by the K\"ahler--Ricci flow \cite{Cao85}, the parabolic complex Monge--Amp\`ere equation \cite{Gil11}, and conformal flows related to the Chern--Yamabe problem \cite{ACS17,LM18,CZ20}. Motivated by these developments, we give a parabolic proof of the semilinear theorem of Fu--Wang--Wu and recover its geometric consequence. On the fixed background $(M,\omega)$, we consider the normalized flow
\begin{equation}\label{eq:intro-flow}
    \partial_t v-\Delta v
    =\psi(|\nabla v|^2)+\langle B,dv\rangle-c(t)-f(x),
\end{equation}
where
\begin{equation*}
    c(t)=
    \frac{\displaystyle\int_M
    \bigl(\Delta v+\psi(|\nabla v|^2)+\langle B,dv\rangle\bigr)
    \,\omega^n}
    {\displaystyle\int_M\omega^n}.
\end{equation*}
Since $f$ has zero mean, this normalization preserves the spatial mean of $v$. Our main results are as follows.

\begin{theorem}[Long-time existence and convergence]\label{thm:semilinear eq}
\label{thm:intro-flow}
Let $(M,\omega)$ be a compact Hermitian manifold, let $B$ be a smooth real $1$-form on $M$, and let $\psi\in C^\infty(\mathbb R)$ satisfy
\begin{equation}\label{eq:psi condition}
    \liminf_{t\to+\infty}\frac{\psi(t)}{t^\mu}\geq\nu>0,
    \qquad \mu>\frac12,
\end{equation}
for constants $\mu$ and $\nu$. For every $v_0\in C^\infty(M)$ and $f\in C^\infty(M)$ with $\int_M f\,\omega^n=0$, the flow \eqref{eq:intro-flow}, with $v(0,\cdot)=v_0$, admits a unique smooth solution for all $t\geq 0$. Moreover, $v(t,\cdot)$ converges uniformly on $M$ to a smooth function $v_\infty$ satisfying
\begin{equation}\label{eq:semilinear eq}
    \Delta v_\infty+\psi(|\nabla v_\infty|^2)
    +\langle B,dv_\infty\rangle=f+c,
\end{equation}
where
\begin{equation*}
    c=c_\infty=
    \frac{\int_M
    \bigl(\Delta v_\infty+\psi(|\nabla v_\infty|^2)
    +\langle B,dv_\infty\rangle\bigr)\,\omega^n}
    {\int_M\omega^n}.
\end{equation*}
The constant $c$ is unique, and the stationary solution is unique up to addition of a constant. The limit is normalized by
\begin{equation*}
    \int_M v_\infty\,\omega^n=\int_M v_0\,\omega^n.
\end{equation*}
\end{theorem}

\begin{theorem}[Geometric consequence of {\cite{FWW13}}]\label{thm:main}
\label{thm:intro-geometric}
Let $(M,\omega)$ be a compact Hermitian manifold of complex dimension $n\geq 2$. For every integer $1\leq k\leq n-1$, there exist a unique constant $\gamma_k$ and a smooth function $v$, unique up to addition of a constant, such that
\begin{equation}\label{eq:main}
    \frac{\sqrt{-1}}{2}\partial\bar\partial(e^v\omega^k)
    \wedge\omega^{n-k-1}
    =\gamma_k e^v\omega^n.
\end{equation}
In particular, when $\gamma_k=0$, the metric $\widetilde\omega=e^{v/k}\omega$ is $k$-Gauduchon.
\end{theorem}

The geometric result follows by taking $\psi(s)=s$, the geometric $1$-form $B$ above, and
\begin{equation*}
    \overline\varphi=
    \frac{\int_M\varphi\,\omega^n}
    {\int_M\omega^n},
    \qquad
    f=\overline\varphi-\varphi,
    \qquad
    \gamma_k(\omega)=\frac{c+\overline\varphi}{n}.
\end{equation*}
Thus the limiting function solves the $\gamma_k$-equation, and the conformal class contains a $k$-Gauduchon metric precisely when $\gamma_k(\omega)=0$.

Our convergence argument does not rely on a gradient-flow structure or a monotone energy functional. The quantity $\partial_t v+c(t)$ satisfies a homogeneous linear uniformly parabolic equation with uniformly bounded drift. The maximum principle and a parabolic Harnack inequality yield a contraction of the spatial oscillation of $\partial_t v$. Its zero mean then implies exponential decay and uniform convergence of $v$. Uniform elliptic $W^{2,p}$ estimates provide strong convergence of the gradients along a subsequence, allowing passage to the nonlinear stationary equation; elliptic regularity then gives a smooth limit.

The paper is organized as follows. Section \ref{sec:short time existence} establishes short-time existence and uniqueness. Section \ref{sec:global existence} derives uniform a priori estimates and obtains global existence by continuation. Section \ref{sec:convergence} proves convergence and identifies the limiting equation. Section \ref{sec:metrics} proves uniqueness and applies the semilinear result to generalized Gauduchon metrics.

\section{Preliminaries and short-time existence}\label{sec:short time existence}
Let $M$ be a compact $n$-dimensional complex manifold. Let $g$ be a Hermitian metric on $M$ and $\omega$ its Hermitian form. 
We write
\begin{equation*}
    \omega = \frac{\sqrt{-1}}{2} \sum_{i,j=1}^{n} g_{i\bar{j}} dz^{i} \wedge d\bar{z}^{j}.
\end{equation*}

Let $(g^{i\bar{j}})$ be the transposed inverse of the matrix $(g_{i\bar{j}})$. 
For any two real 1-forms $A$ and $B$ on $M$, locally given by
\begin{equation*}
    A = \sum_{i=1}^{n} (A_{i} dz_{i} + A_{\bar{i}} d\bar{z}_{i})
    \quad \text{and} \quad
    B = \sum_{i=1}^{n} (B_{i} dz_{i} + B_{\bar{i}} d\bar{z}_{i}),
\end{equation*}
we denote
\begin{equation*}
    \langle A, B \rangle= \frac{1}{2} \sum_{i,j=1}^{n} g^{i\bar{j}} (A_{i} B_{\bar{j}} + A_{\bar{j}} B_{i}).
\end{equation*}
In particular, we have
\begin{equation*}
    \langle dh, dh \rangle= \sum_{i,j=1}^{n} g^{i\bar{j}} \frac{\partial h}{\partial z_{i}} \frac{\partial h}{\partial \bar{z}_{j}}
    \equiv |\nabla h|^{2}
    \quad \text{for all } h \in C^{1}(M).
\end{equation*}
The Laplacian $\Delta$ associated with $\omega$ is given by
\begin{equation*}
    \Delta h
    = \frac{n \omega^{n-1} \wedge (\sqrt{-1}/2) \partial \bar{\partial} h}{\omega^{n}}
    = \sum_{i,j=1}^{n} g^{i\bar{j}} h_{i\bar{j}}
    \quad \text{for all } h \in C^{2}(M).
\end{equation*}

In the general Hermitian setting, the real Laplacian $\Delta_g$ and complex Laplacian $\Delta$ defined here differ and satisfy the following identity
\begin{equation}\label{eq:difference of real and complex laplacian}
    \Delta_g v=4(\Delta v+\langle B_1,dv\rangle),
\end{equation}
where 
\begin{equation*}
    B_1=\frac{\sqrt{-1}}{(n-1)!}*(\partial\omega^{n-1}-\bar\partial\omega^{n-1}).
\end{equation*}
This identity, together with the well-developed theory of the real Laplacian, underlies both the short-time existence proof in the remainder of this section and the convergence argument in Section \ref{sec:convergence}. 

In order to prove Theorem \ref{thm:semilinear eq}, we consider the following heat flow:
\begin{equation}\label{eq:heat flow v}
    \left\{\begin{gathered}
        \partial_t v-\Delta v=\psi(\vert\nabla v\vert^2)+\langle B,dv\rangle-c(t)-f(x), \quad x\in M,\,t\ge 0, \\
        v(0,x)= v_0(x) , \quad x\in M.
\end{gathered}\right.
\end{equation}
Here $\psi$, $f$ and $B$ are the same as in Theorem \ref{thm:semilinear eq}, $v_0\in C^{\infty}(M)$ and
\begin{equation*}
    \begin{split}
        c(t)
        =&\frac{\int_{M}(\Delta v+\psi(\vert\nabla v\vert^2)+\langle B,dv\rangle)\omega^n}{\int_{M}\omega^n} \\
        =&\frac{\int_{M}(\psi(\vert\nabla v\vert^2)+\langle B-B_1,dv\rangle)\omega^n}{\int_{M}\omega^n}
    \end{split}
\end{equation*}
Note that $c=c(t)$  depends only on $t$, i.e., it is spatially constant for each fixed $t$. In this setting,
\begin{equation*}
    \int_{M} \partial_t v(t)\ \omega^n\equiv0, \quad\forall t\ge 0.
\end{equation*}
Since adding a constant to $v$ does not change the equation, we may assume without loss of generality that
\begin{equation*}
    \int_{M} v(t)\ \omega^n\equiv0, \quad\forall t\ge 0.
\end{equation*}

\begin{prop}\label{thm:short time existence}
Let $(M,\omega)$ be a compact Hermitian manifold. The heat flow (\ref{eq:heat flow v}) admits a unique solution 
\begin{equation*}
    v\in C([0,T],C^1(M))\cap C^{\infty}((0,T]\times M).
\end{equation*}
Here $T>0$ admits a positive lower bound depending only on $\|v_0\|_{C^1(M)}$ and the fixed data $M,\omega,B,\psi,f$. In particular, a common positive existence time can be chosen when the initial data range over a bounded subset of $C^1(M)$.
\end{prop}

\begin{proof}
Set $t^{\prime}=\frac{t}{4}$, and $u(t^{\prime},\cdot)=v(4t^{\prime},\cdot)$. $v(t,x)$ solves (\ref{eq:heat flow v}) on $[0,T]\times M$ if and only if $u(t^{\prime},x)$ solves the following system
\begin{equation}\label{eq:heat equation}
    \left\{\begin{gathered}
                 (\frac{\partial}{\partial t^{\prime}}-\Delta_g) u =4F(x,t^{\prime},\nabla u), \quad \text{on}\ [0,\frac{T}{4}]\times M\\
                 u(0,x)=v_0(x), \quad x\in M.
           \end{gathered}\right.
\end{equation}
where $F(x,t^{\prime},\nabla u)=\psi(\vert\nabla u\vert^2)+\langle B^{\prime},du\rangle-c(4t^{\prime})-f(x)$ and $B^{\prime}=B-B_1$.
Denote $\Phi(u)=4F(x,t^{\prime},\nabla u)$.
We introduce the following operator:
\begin{equation}\label{eq:solution operator}
    \begin{split}
        \Psi(u(x,t^{\prime}))
        &=e^{t^{\prime}\Delta_g} v_0+\int^{t^{\prime}}_0 e^{(t^{\prime}-s)\Delta_g}\Phi(u(s))ds \\
        &=\int_{M} K(x,y,t^{\prime})v_0(y)\,\mathrm{d}y+\int^{t^{\prime}}_0\int_{M} K(x,y,t^{\prime}-s)\Phi(u(s))(y)\,\mathrm{d}y\,\mathrm{d}s,
    \end{split}
\end{equation}
where $\{e^{t^{\prime}\Delta_g},t^{\prime}>0\}$ is the heat semigroup and $K(x,y,t^{\prime})$ is the heat kernel for $\Delta_g$ on $M$. 
It is clear that the fixed point of $\Psi$ solves (\ref{eq:heat equation}). To complete the proof, we only need to show that $\Psi$ has a unique fixed point $u$. The heat semigroup is strongly continuous on $C^1(M)$.
By the uniform boundedness principle and the standard heat
semigroup smoothing estimate, there exists a constant
$C_1\ge1$, depending only on the background geometry, such that
\begin{equation}\label{eq:short-time-C1-semigroup}
\begin{aligned}
\sup_{0\le t'\le1}
\|e^{t'\Delta_g}\|_{\mathcal L(C^1(M),C^1(M))}
&\le C_1,\\
\|e^{t'\Delta_g}\|_{\mathcal L(C(M),C^1(M))}
&\le C_1(t')^{-1/2},
\qquad 0<t'\le1.
\end{aligned}
\end{equation}
For $0<T\le4$, set
\begin{equation*}
\begin{split}
Z=\biggl\{
u\in C([0,T/4],C^1(M)):\;&u(0)=v_0,\\
&\sup_{0\le t'\le T/4}
\|u(t')-e^{t'\Delta_g}v_0\|_{C^1(M)}\le\alpha
\biggr\},
\end{split}
\end{equation*}
where $\alpha>0$ is fixed. This is a nonempty closed subset of $C([0,T/4],C^1(M))$, equipped with the supremum norm.
For every $u\in Z$ and $0\le t'\le T/4$, we have
\begin{equation*}
\|u(t')\|_{C^1(M)}
\le \alpha+C_1\|v_0\|_{C^1(M)}.
\end{equation*}

Using a $C^1$ norm that controls both the function and its gradient, we may take
\begin{equation*}
C_2=
\max_{0\le s\le(\alpha+C_1\|v_0\|_{C^1(M)})^2}
|\psi(s)|.
\end{equation*}
Recall that $B'=B-B_1$ and that the spatially constant term in $\Phi$ is the average of $\psi(|\nabla u|^2)+\langle B',du\rangle$.
It follows that
\begin{equation}\label{eq:short-time-nonlinearity-bound}
\begin{aligned}
\|\Phi(u(t'))\|_{C(M)}
&\le 8C_2+
8\max_{x\in M}|B'|
\bigl(\alpha+C_1\|v_0\|_{C^1(M)}\bigr)
+4\max_{x\in M}|f|\\
&=:C_3.
\end{aligned}
\end{equation}
Set
\begin{equation*}
\begin{aligned}
C_4={}&16\bigl(\alpha+C_1\|v_0\|_{C^1(M)}\bigr)
\max_{0\le s\le(\alpha+C_1\|v_0\|_{C^1(M)})^2}
|\psi'(s)|\\
&+8\max_{x\in M}|B'|.
\end{aligned}
\end{equation*}
Using the smoothness of $\psi$ and estimating the difference of the spatial averages by the supremum norm of the difference of their integrands, we obtain, for $u_1,u_2\in Z$,
\begin{equation}\label{eq:short-time-lipschitz}
\|\Phi(u_1(t'))-\Phi(u_2(t'))\|_{C(M)}
\le C_4\|u_1(t')-u_2(t')\|_{C^1(M)}.
\end{equation}

The Duhamel term is continuous with values in $C^1(M)$, including at $t'=0$.
By \eqref{eq:short-time-C1-semigroup} and
\eqref{eq:short-time-nonlinearity-bound}, for $u\in Z$,
\begin{equation}\label{eq:short-time-self-map}
\begin{aligned}
\|\Psi(u)(t')-e^{t'\Delta_g}v_0\|_{C^1(M)}
&=\left\|\int_0^{t'}e^{(t'-s)\Delta_g}
           \Phi(u(s))\,ds\right\|_{C^1(M)}\\
&\le 2C_1C_3(t')^{1/2}.
\end{aligned}
\end{equation}
Similarly, by \eqref{eq:short-time-C1-semigroup} and \eqref{eq:short-time-lipschitz},
\begin{equation}\label{eq:short-time-contraction}
\begin{aligned}
&\|\Psi(u_1)-\Psi(u_2)\|_{C([0,T/4],C^1(M))}\\
&\qquad\le
2C_1C_4(T/4)^{1/2}
\|u_1-u_2\|_{C([0,T/4],C^1(M))}.
\end{aligned}
\end{equation}
Choose
\begin{equation*}
T=4\min\left\{
1,\,
\frac{\alpha^2}{16C_1^2(1+C_3)^2},\,
\frac{1}{16C_1^2(1+C_4)^2}
\right\}.
\end{equation*}
Then
\begin{equation*}
2C_1C_3(T/4)^{1/2}\le\frac{\alpha}{2},
\qquad
2C_1C_4(T/4)^{1/2}\le\frac12.
\end{equation*}
Thus $\Psi$ maps $Z$ into itself and is a contraction. By the contraction mapping theorem, $\Psi$ has a unique fixed point $u\in Z$, which is a mild solution of (\ref{eq:heat equation}).

Uniqueness in the full class $C([0,T/4],C^1(M))$ of mild solutions follows by applying the same difference estimate on successive sufficiently short subintervals, using a common $C^1$ bound for any two solutions under consideration.

Returning to the original time variable and using parabolic smoothing and bootstrapping, we obtain
\begin{equation*}
v\in C([0,T],C^1(M))\cap C^\infty((0,T]\times M).
\end{equation*}
All the constants used to choose $T$ depend only on $\alpha,\|v_0\|_{C^1(M)},M,\omega,B,\psi,f$. Fixing, for example, $\alpha=1$, we obtain a positive existence time whose dependence on the initial datum is only through its $C^1(M)$ norm. The displayed choices of the constants also show that a common positive existence time works for initial data in any bounded subset of $C^1(M)$.
\end{proof}

\section{A priori estimates and global existence}\label{sec:global existence}

First, we intend to give an estimate for $\partial_{t} v$.
\begin{lemma}\label{thm:estimate of t derivative}
    Let $v\in C^3([0,T)\times M)$ be a solution of (\ref{eq:heat flow v}) on $[0,T)\times M$. For $0\le t<t+s< T$, we have
    \begin{equation}\label{eq:min non-decrease}
        \underset{x\in M}{\text{min}}(\partial_t v+c)(t)\le\underset{x\in M}{\text{min}}(\partial_t v+c)(t+s),
    \end{equation}
    and
    \begin{equation}\label{eq:max non-increase}
        \underset{x\in M}{\text{max}}(\partial_t v+c)(t+s)\le\underset{x\in M}{\text{max}}(\partial_t v+c)(t).
    \end{equation}
     Moreover, $\operatorname*{osc}_{x\in M}(\partial_t v)$ is monotone non-increasing with respect to $t$, and
    \begin{equation*}
        \underset{(t,x)\in[0,T)\times M}{\text{sup}}\vert\partial_{t} v\vert\le C,
    \end{equation*}
    where $C>0$ is a constant depending only on $B$, $f$, $\omega$, $\psi$ and $v_0$.
\end{lemma}
\begin{proof}
Differentiating (\ref{eq:heat flow v}) with respect to t, we get
\begin{equation}\label{eq:partial_t of semileaner eq}
    \partial_t (\partial_t v)-\Delta (\partial_t v)=2\psi^{\prime}(\vert\nabla v\vert^2)\langle dv,d(\partial_t v)\rangle+\langle B,d(\partial_t v)\rangle-\partial_tc.   
\end{equation}
Note that $c$ only depends on $t$. In fact, we have
\begin{equation}\label{eq:t derivative of the flow}
    \partial_t (\partial_t v+c)-\Delta (\partial_t v+c)=2\psi^{\prime}(\vert\nabla v\vert^2)\langle dv,d(\partial_t v+c)\rangle+\langle B,d(\partial_t v+c)\rangle. 
\end{equation}
Applying maximum principle to $\partial_t v+c$ in (\ref{eq:t derivative of the flow}) implies
\begin{equation}\label{eq:later time control}
    \underset{x\in M}{\text{min}}(\partial_t v+c)(t)\le(\partial_t v+c)(t+s)\le\underset{x\in M}{\text{max}}(\partial_t v+c)(t).
\end{equation}
It leads to
    \begin{equation*}
        \underset{x\in M}{\text{min}}(\partial_t v+c)(t)\le\underset{x\in M}{\text{min}}(\partial_t v+c)(t+s),
    \end{equation*}
    and
    \begin{equation*}
        \underset{x\in M}{\text{max}}(\partial_t v+c)(t+s)\le\underset{x\in M}{\text{max}}(\partial_t v+c)(t).
    \end{equation*}
Subtracting the two inequalities above
    \begin{equation*}
        \operatorname*{osc}_{x\in M}(\partial_t v)(t+s)\le\operatorname*{osc}_{x\in M}(\partial_t v)(t).
    \end{equation*}
Since $\int_M\partial_t v\,\omega^n=0$, we have
\begin{equation*}
c(t)=
\frac{\int_M(\partial_t v+c)(t,x)\,\omega^n}
     {\int_M\omega^n}.
\end{equation*}
Applying \eqref{eq:later time control} with initial time zero and using this identity, we obtain
\begin{equation}\label{eq:uniform bounds for time derivative and c}
\begin{aligned}
|c(t)|
&\leq \|(\partial_t v+c)(0,\cdot)\|_{C^0(M)},\\
\|\partial_t v(t,\cdot)\|_{C^0(M)}
&\leq \operatorname{osc}_M(\partial_t v+c)(0,\cdot).
\end{aligned}
\end{equation}
Since
\begin{equation*}
(\partial_t v+c)(0,\cdot)
=\Delta v_0+\psi(|\nabla v_0|^2)+\langle B,dv_0\rangle-f,
\end{equation*}
both bounds depend only on $B$, $f$, $\omega$, $\psi$ and $v_0$. This completes the proof.
\end{proof}
Next, we derive the $C^1$ estimate for $v$. More precisely, we have the following Lemma.
\begin{lemma}\label{thm:C_1 estimate}
    Let $v\in C^3([0,T)\times M)$ be a solution of (\ref{eq:heat flow v}) on $[0,T)\times M$. We have
    \begin{equation*}
        \underset{(t,x)\in [0,T)\times M}{\text{sup}}\vert\nabla v\vert\le C,
    \end{equation*}
    where $C>0$ is a constant depending only on $B$, $f$, $\omega$, $\psi$, $\mu$, $\nu$ and $v_0$.
\end{lemma}
\begin{proof}
    Throughout the proof, $C>0$ denotes a generic uniform constant depending only on $B$, $f$, $\omega$, $\psi$, $\mu$, $\nu$ and $v_0$, which may vary from line to line.
    Consider the linear parabolic operator
    \begin{equation*}
        L(h)=-\partial_t h+\Delta h+2\psi^{\prime}(\vert \nabla v\vert^2)\langle dh,dv\rangle_{\omega}.
    \end{equation*}
    Fix $0<T'<T$. Since $[0,T']\times M$ is compact,
$|\nabla v|^2$ attains its maximum at some point
$(t_0,x_0)\in[0,T']\times M$.
If $t_0=0$, then
\begin{equation*}
    \sup_{(t,x)\in[0,T']\times M}|\nabla v(t,x)|^2
    =|\nabla v_0(x_0)|^2
    \leq \|\nabla v_0\|_{C^0(M)}^2,
\end{equation*}
which gives the desired estimate. It therefore remains to consider the case $0<t_0\leq T'$. At $(t_0,x_0)$, the spatial maximum conditions give
\begin{equation*}
    \Delta(|\nabla v|^2)(t_0,x_0)\leq 0,
    \qquad
    \nabla(|\nabla v|^2)(t_0,x_0)=0.
\end{equation*}
Moreover,
\begin{equation*}
    \partial_t(|\nabla v|^2)(t_0,x_0)\geq 0.
\end{equation*}
Indeed, this derivative vanishes if $0<t_0<T'$, while it is nonnegative by taking the derivative from the left if $t_0=T'$. Consequently, we have
    \begin{equation}\label{eq:maximum point}
        L(\vert \nabla v\vert^2)(t_0,x_0)\le 0.
    \end{equation}
    
    Straightforward computations give that 
    \begin{equation*}
        \begin{split}
            L(\vert\nabla v\vert^2)=&-\partial_t(\vert\nabla v\vert^2)+\Delta(\vert\nabla v\vert^2)+\psi^{\prime}g^{i\bar{j}}[(\vert\nabla v\vert^2)_{i}v_{\bar{j}}+v_{i}(\vert\nabla v\vert^2)_{\bar{j}}] \\
            =& -g^{p\bar{q}}(v_{tp}v_{\bar{q}}+v_{p}v_{t\bar{q}})+g^{i\bar{j}}g^{p\bar{q}}(v_{pi}v_{\bar{q}\bar{j}}+v_{p\bar{j}}v_{i\bar{q}})+g^{p\bar{q}}[(\Delta v)_{p}v_{\bar{q}}+v_{p}(\Delta v)_{\bar{q}}]\\
            &+g^{i\bar{j}}g^{p\bar{q}}_{,i}(v_{p\bar{j}}v_{\bar{q}}+v_{p}v_{\bar{q}\bar{j}})+g^{i\bar{j}}g^{p\bar{q}}_{,\bar{j}}(v_{pi}v_{\bar{q}}+v_{p}v_{i\bar{q}})\\
            &-g^{p\bar{q}}(g^{i\bar{j}}_{,p}v_{i\bar{j}}v_{\bar{q}}+g^{i\bar{j}}_{,\bar{q}}v_{p}v_{i\bar{j}})+\psi^{\prime}g^{i\bar{j}}[(\vert\nabla v\vert^2)_{i}v_{\bar{j}}+v_{i}(\vert\nabla v\vert^2)_{\bar{j}}] \\
            &+g^{i\bar{j}}g^{p\bar{q}}_{,i\bar{j}}(v_{p}v_{\bar{q}})
        \end{split}
    \end{equation*}
    Without loss of generality, we take a normal coordinate system around $x_0$ such that
    \begin{equation*}
        g_{i\bar{j}}(x_0)=\delta_{ij} \quad \forall \ 1\le i,j\le n.
    \end{equation*}
    By virtue of Schwarz inequality, we have
    \begin{equation*}
        \begin{split}
             L(\vert\nabla v\vert^2)\ge &\frac{1}{2}\sum_{i,p=1}^{n}\vert v_{p\bar{i}}\vert^2-C\vert\nabla v\vert^2-C \\
             \ge & \frac{1}{2}\sum_{i=1}^{n}\vert v_{i\bar{i}}\vert^2-C\vert\nabla v\vert^2-C \\
             \ge & \frac{1}{2n}\vert\Delta v\vert^2- C\vert\nabla v\vert^2-C\\
             \ge & \frac{1}{2n} \vert\psi(\vert \nabla v\vert^2)+\langle B, dv\rangle-f-c-\partial_t v\vert^2- C\vert\nabla v\vert^2-C\\
             \ge & \frac{1}{4n} \vert\psi(\vert \nabla v\vert^2)\vert^2-C\vert\nabla v\vert^2-C(1+\vert c\vert^2).
        \end{split}  
    \end{equation*}
    We may assume that $\vert\nabla v\vert^2(t_0,x_0)$ is sufficiently large so that at point $(t_0,x_0)$
    \begin{equation*}
        \psi(\vert\nabla v\vert^2)\ge \frac{\nu}{2}\vert\nabla v\vert^{2\mu},
    \end{equation*}
    where $\mu>\frac{1}{2}$ and $\nu>0$ are constants in (\ref{eq:psi condition}). Combining (\ref{eq:maximum point}) with the preceding estimates and the uniform bound for $c$ in Lemma \ref{thm:estimate of t derivative}, we obtain
\begin{equation*}
    |\nabla v|^2(t_0,x_0)\leq C
\end{equation*}
when $t_0>0$. Here the case in which $|\nabla v|^2(t_0,x_0)$ is below the threshold required for the growth estimate of $\psi$ is included by enlarging $C$.
Together with the case $t_0=0$, this yields
\begin{equation*}
    \sup_{(t,x)\in[0,T']\times M}|\nabla v(t,x)|^2
    \leq \max\left\{
        \|\nabla v_0\|_{C^0(M)}^2,C
    \right\}.
\end{equation*}
Since the right-hand side is independent of $T'$, letting $T'\uparrow T$ proves the lemma.
\end{proof}

The $C^0$ estimate is then a direct corollary of the $C^1$ estimate.
\begin{corollary}\label{thm:C_0 estimate}
        Let $v\in C^3([0,T)\times M)$ be a solution of (\ref{eq:heat flow v}) on $[0,T)\times M$. We have
    \begin{equation*}
        \underset{(t,x)\in [0,T)\times M}{\text{sup}}\vert v\vert\le C,
    \end{equation*}
    where $C>0$ is a constant depending only on $B$, $f$, $\omega$, $\psi$, $\mu$, $\nu$, $v_0$ and $M$.
\end{corollary}

\begin{proof}
    Since
    \begin{equation*}
        \int_M  v=0, \ \forall t\in [0,T),
    \end{equation*}
    we know for any $t\in [0,T)$ that there exists some $x^{\prime}\in M$ such that $v(t,x^{\prime})=0$. Using the Hopf-Rinow Theorem, for any $x\in M$, there is a minimizing geodesic curve $\gamma$ connecting $x^{\prime}$ to $x$. Note that the length of the geodesic $\vert\gamma\vert \le diam_{\omega}(M)$. Therefore, for any $(t,x)\in [0,T)\times M$, 
    \begin{equation*}
        \begin{split}
            \vert v(t,x)\vert&=\vert v(t,x)-v(t,x^{\prime})\vert= \left\vert\int_{0}^{1}\frac{d(v\circ\gamma)}{ds}\right\vert \\&\le diam_{\omega}(M)  \int_{0}^{1}(\vert\nabla v\vert\circ\gamma) ds \le C.
        \end{split}
    \end{equation*}
\end{proof}
    Using Proposition \ref{thm:short time existence}, Lemma \ref{thm:C_1 estimate} and Corollary \ref{thm:C_0 estimate}, we have 
    \begin{prop}\label{thm:long time existence}
        There exists a unique function $v\in C([0,+\infty),C^1(M))\cap C^{\infty}((0,+\infty)\times M)$ solving (\ref{eq:heat flow v}).
    \end{prop}

\section{Convergence of the flow}\label{sec:convergence}

Our proof of the convergence of the flow (\ref{eq:heat flow v}) is based on a global Harnack inequality. To this end, we first recall a corollary of Saloff-Coste \cite{Sal92}, which plays a key role in our argument. 
\begin{lemma}\label{thm:local harnack}\cite[Corollary 5.4]{Sal92}
    Suppose $M$ is an n-dimensional complete Riemannian manifold. Let $A_t$ be a measurable section of $\operatorname{End}(TM)$ depending on $t$ and let $m$ be a positive measurable function on $M$. Suppose 
    \begin{equation*}
        \mu^{-1} \leq m(x) \leq \mu, \qquad x \in M, 
    \end{equation*}
    and
    \begin{equation*}
        \alpha^{-1}|X|^{2} \leq g(A_{t,x}X, X), \quad |A_{t,x}X| \leq \alpha|X| \quad \forall t \in \mathbb{R}, x \in M, X \in TM_{x}.
    \end{equation*}
    Here $\mu$ and $\alpha$ are positive constants.
    
    Let $X_t$ and $Y_t$ be two vector fields on $M$ such that the functions $(t,x)\to (x,X_{t,x})$ and $(t,x)\to (x,Y_{t,x})$ are measurable, and let $b(t,x)$ be a measurable function on $\mathbb{R}\times M$. We assume that there exists a constant $\beta\ge 0$ such that
    \begin{equation*}
        |X|^{2} + |Y|^{2} + |b| \leq \beta.
    \end{equation*}

    If $u$ is a positive solution of the following equation 
    \begin{equation}\label{eq:saloff_coste}
        \partial_t u -m^{-1} \operatorname{div}\left(mA( \nabla u + u X)\right) -\langle Y,\nabla u\rangle - b u=0 \quad \text{on}\ (0,T)\times B(x,r),
    \end{equation}
    then, for every $\delta\in (0,1)$, $0<t^{\prime}<t<T$, any $y,y^{\prime}\in B(x,\delta r)\subset B(x,r)$, it satisfies
    \begin{equation*}
        \ln\left(\frac{u(t', y')}{u(t, y)}\right) \leq C\left(\frac{\rho^{2}}{(t - t')} + \left(\beta + K + \frac{1} {r^{2}} + \frac{1}{t'}\right)(t - t')+1\right).
    \end{equation*}
    Here $\rho=\rho(y,y^{\prime})$, $-K$ is the lower bound of Ricci curvature in $B(x,2r)$, and $C$ is a constant only depending on $n$, $\delta$, $\alpha$, and $\mu$.
\end{lemma}

   \begin{remark}
     The printed statement of \cite[Corollary~5.4]{Sal92} omits the final additive term \(+1\) inside the brackets on the right-hand side. This omission is explicitly pointed out by Colombo--Gama--Mari--Rigoli in \cite[Appendix~A, formula~(69) and the discussion following it]{CGMR24}. We use the corrected estimate throughout. In the Harnack-chain argument below, the contribution of this additional term is absorbed into the constant \(C_H\).
   \end{remark}

    In what follows, we regard $M$ as a real $2n$-dimensional Riemannian manifold equipped with the metric associated with $\omega$, and apply Lemma~\ref{thm:local harnack} with dimension $2n$. For the equation under consideration, we take $m\equiv1$, $A=\frac14\mathrm{Id}$, $X\equiv0$, and $b\equiv0$, with $\alpha=4$ and $\mu=1$. The first-order term is given by
    \begin{equation*}
        \langle Y,\nabla u\rangle=2\psi'(|\nabla v|^2)\langle dv,du\rangle+\langle B-B_1,du\rangle,
    \end{equation*}
where the left-hand side uses the real Riemannian gradient and pairing, while the pairings on the right-hand side are those defined in Section \ref{sec:short time existence}. Thus \eqref{eq:saloff_coste} takes the form
    \begin{equation}\label{eq:t derivative of the flow 2}
        \partial_t u-\frac14\Delta_g u=2\psi'(|\nabla v|^2)\langle dv,du\rangle+\langle B-B_1,du\rangle.
    \end{equation}
Moreover, the uniform gradient estimate for $v$ in Section \ref{sec:global existence}, together with the smoothness of $\psi$, yields
    \begin{equation*}
         |Y|^2\leq\beta\qquad\text{on }[0,+\infty)\times M,
    \end{equation*}
for some constant $\beta\geq0$ depending only on the fixed background data and the initial datum. The required bound on the lower-order coefficients therefore follows.
 
\begin{lemma}\label{thm:global harnack}
Let $T_0\geq 0$, and fix $\sigma,\tau>0$.
Suppose that $u>0$ solves (\ref{eq:t derivative of the flow 2}) on $(T_0,\infty)\times M$, and that the coefficient bounds above hold uniformly in time. Then, for every $t\geq T_0+\sigma$,
\begin{equation*}
    \max_{x\in M}u(t,x)
    \leq C_H\min_{x\in M}u(t+\tau,x),
\end{equation*}
where
$C_H=C_H(n,\delta,K,\beta,\operatorname{diam}_{\omega}(M),
\sigma,\tau)$ is independent of $T_0$, $t$, and $u$.
\end{lemma}

\begin{proof}
Set $d=\operatorname{diam}_{\omega}(M)$ and $r=d/2$.
Fix $\delta\in(0,1)$, and choose
\begin{equation*}
    N=\left\lfloor\frac{4}{\delta}\right\rfloor+1.
\end{equation*}
Here $\left\lfloor\cdot\right\rfloor$ is the floor function. For arbitrary $y,y'\in M$, let $\gamma:[0,\ell]\to M$ be a unit-speed minimizing geodesic joining $y$ to $y'$, where $\ell\leq d$. 
Define
\begin{equation*}
    y_i=\gamma\left(\frac{i\ell}{N}\right),
    \qquad
    t_i=t+\frac{i\tau}{N},
    \qquad i=0,\ldots,N.
\end{equation*}
If $y=y'$, take $y_i=y$ for all $i$.
Then
\begin{equation*}
    \operatorname{dist}(y_{i-1},y_i)
    =\frac{\ell}{N}<\delta r.
\end{equation*}

Apply the local Harnack inequality to the time-translated solution
\begin{equation*}
    \widetilde u(s,x)=u(T_0+s,x),
    \qquad s>0,
\end{equation*}
on each ball $B(y_{i-1},r)$. Since the translated coefficients satisfy the same uniform bounds, we obtain
\begin{equation*}
\begin{aligned}
    \log\frac{u(t_{i-1},y_{i-1})}{u(t_i,y_i)}
    &\leq C_0\left[
        1+\frac{(\ell/N)^2}{\tau/N}
        +\left(
            \beta+K+\frac{1}{r^2}
            +\frac{1}{t_{i-1}-T_0}
        \right)\frac{\tau}{N}
    \right] \\
    &\leq C_0\left[
        1+\frac{\ell^2}{N\tau}
        +\left(
            \beta+K+\frac{4}{d^2}+\frac{1}{\sigma}
        \right)\frac{\tau}{N}
    \right],
\end{aligned}
\end{equation*}
where we used $t_{i-1}-T_0\geq\sigma$. The constant $C_0$ depends only on the parameters in the local Harnack inequality.

Summing these inequalities over $i=1,\ldots,N$ gives
\begin{equation*}
    \log\frac{u(t,y)}{u(t+\tau,y')}
    \leq C_0\left[
        N+\frac{d^2}{\tau}
        +\left(
            \beta+K+\frac{4}{d^2}+\frac{1}{\sigma}
        \right)\tau
    \right].
\end{equation*}
Thus we may take
\begin{equation*}
    C_H=
    \exp\left\{
        C_0\left[
            N+\frac{d^2}{\tau}
            +\left(
                \beta+K+\frac{4}{d^2}+\frac{1}{\sigma}
            \right)\tau
        \right]
    \right\}.
\end{equation*}
Taking the maximum over $y$ and the minimum over $y'$ completes the proof.
\end{proof}

With the global Harnack inequality at hand, we can derive a decay estimate for the oscillation of $\partial_t v$.
\begin{theorem}\label{thm:time derivative control}
    Fix $\tau>0$. There exists a constant $0<\theta(C_H)<1$ independent of $t$, such that for every $t>0$, the following inequality holds
    \begin{equation}\label{eq:time derivative osc control}
        \operatorname*{osc}_{x\in M}(\partial_t v(t+\tau,\cdot))\le \theta \ \operatorname*{osc}_{x\in M}(\partial_t v(t,\cdot)).
    \end{equation}
    Here $\operatorname*{osc}_{x\in M}(u)=\underset{x\in M}{\text{max}}(u)-\underset{x\in M}{\text{min}}(u)$. 
    Moreover, we have
    \begin{equation}\label{eq:time derivative sup control}
        \sup_{x\in M} |\partial_t v(t,\cdot)| \leq \operatorname*{osc}_{x\in M}\left(\Delta v_{0} + \psi(|\nabla v_{0}|^{2}) + \langle B, dv_{0} \rangle - f\right) \theta^{\frac{t}{\tau} -2}
    \end{equation}
    holds for every $t\ge 2\tau$.
\end{theorem}

\begin{proof}
We first prove the oscillation contraction for an arbitrary starting time $T_0>0$. Denote
\begin{equation*}
M_0:=\max_{x\in M}\partial_t v(T_0,\cdot)+c(T_0),
\qquad
m_0:=\min_{x\in M}\partial_t v(T_0,\cdot)+c(T_0).
\end{equation*}
For $t\geq T_0$, define
\begin{equation*}
w_1(t):=M_0-(\partial_t v(t,\cdot)+c(t))+\varepsilon,
\qquad
w_2(t):=(\partial_t v(t,\cdot)+c(t))-m_0+\varepsilon,
\end{equation*}
where $\varepsilon>0$ is a constant.
By Lemma \ref{thm:estimate of t derivative}, we have
\begin{equation*}
w_1(t),w_2(t)\geq\varepsilon>0,
\qquad t\geq T_0.
\end{equation*}
Moreover, by (\ref{eq:t derivative of the flow}), $w_1$ and $w_2$ satisfy
\begin{equation*}
\partial_t w_1-\Delta w_1
=
\left\langle
2\psi'(|\nabla v|^2)dv+B,dw_1
\right\rangle
\end{equation*}
and
\begin{equation*}
\partial_t w_2-\Delta w_2
=
\left\langle
2\psi'(|\nabla v|^2)dv+B,dw_2
\right\rangle,
\end{equation*}
respectively.

Applying the Harnack inequality as in Lemma \ref{thm:global harnack}, after translating the initial time $T_0$ to zero, we obtain
\begin{equation*}
\max_{x\in M}w_1(T_0+\tau/2,\cdot)
\leq
C_H\min_{x\in M}w_1(T_0+\tau,\cdot)
\end{equation*}
and
\begin{equation*}
\max_{x\in M}w_2(T_0+\tau/2,\cdot)
\leq
C_H\min_{x\in M}w_2(T_0+\tau,\cdot).
\end{equation*}
Here the earlier comparison time is a fixed positive distance $\tau/2$ from the initial time. Since the coefficient bounds are uniform in time, $C_H$ depends only on $\tau$, the background geometry, and these bounds, and is independent of $T_0$ and $\varepsilon$. Enlarging $C_H$ if necessary, we assume $C_H>1$.

Note that
\begin{equation*}
\begin{aligned}
\max_{x\in M}w_1(T_0+\tau/2,\cdot)
&=
M_0-\min_{x\in M}\partial_t v(T_0+\tau/2,\cdot)
-c(T_0+\tau/2)+\varepsilon,\\
\min_{x\in M}w_1(T_0+\tau,\cdot)
&=
M_0-\max_{x\in M}\partial_t v(T_0+\tau,\cdot)
-c(T_0+\tau)+\varepsilon.
\end{aligned}
\end{equation*}
Thus,
\begin{equation*}
\begin{aligned}
&M_0-\min_{x\in M}\partial_t v(T_0+\tau/2,\cdot)
-c(T_0+\tau/2)+\varepsilon\\
&\qquad\leq
C_H\Bigl(
M_0-\max_{x\in M}\partial_t v(T_0+\tau,\cdot)
-c(T_0+\tau)+\varepsilon
\Bigr).
\end{aligned}
\end{equation*}
Similarly, applying the Harnack inequality to $w_2$ gives
\begin{equation*}
\begin{aligned}
&\max_{x\in M}\partial_t v(T_0+\tau/2,\cdot)
+c(T_0+\tau/2)-m_0+\varepsilon\\
&\qquad\leq
C_H\Bigl(
\min_{x\in M}\partial_t v(T_0+\tau,\cdot)
+c(T_0+\tau)-m_0+\varepsilon
\Bigr).
\end{aligned}
\end{equation*}
Adding these two inequalities, we obtain
\begin{equation*}
\begin{aligned}
&M_0-m_0
+\operatorname{osc}_{x\in M}
  \bigl(\partial_t v(T_0+\tau/2,\cdot)\bigr)
+2\varepsilon\\
&\qquad\leq
C_H\Bigl(
M_0-m_0
-\operatorname{osc}_{x\in M}
  \bigl(\partial_t v(T_0+\tau,\cdot)\bigr)
+2\varepsilon
\Bigr).
\end{aligned}
\end{equation*}
Letting $\varepsilon\downarrow0$ and using the monotonicity of $\operatorname{osc}_{x\in M}(\partial_t v)$, we find
\begin{equation*}
(C_H+1)
\operatorname{osc}_{x\in M}
\bigl(\partial_t v(T_0+\tau,\cdot)\bigr)
\leq
(C_H-1)(M_0-m_0).
\end{equation*}
Since
\begin{equation*}
M_0-m_0
=\operatorname{osc}_{x\in M}
\bigl(\partial_t v(T_0,\cdot)\bigr),
\end{equation*}
it follows that
\begin{equation*}
\operatorname{osc}_{x\in M}
\bigl(\partial_t v(T_0+\tau,\cdot)\bigr)
\leq
\theta\,
\operatorname{osc}_{x\in M}
\bigl(\partial_t v(T_0,\cdot)\bigr),
\qquad
\theta:=\frac{C_H-1}{C_H+1}\in(0,1).
\end{equation*}
As $T_0>0$ was arbitrary and $\theta$ is independent of $T_0$, we obtain
\begin{equation*}
\operatorname{osc}_{x\in M}
\bigl(\partial_t v(t+\tau,\cdot)\bigr)
\leq
\theta\,
\operatorname{osc}_{x\in M}
\bigl(\partial_t v(t,\cdot)\bigr),
\qquad t>0.
\end{equation*}
This proves (\ref{eq:time derivative osc control}).

For $t\geq2\tau$, iterating the preceding inequality $\lfloor t/\tau\rfloor-1$ times yields
\begin{equation*}
\begin{aligned}
&\operatorname{osc}_{x\in M}
\bigl(\partial_t v(t,\cdot)\bigr)\\
&\quad\leq
\theta^{\lfloor t/\tau\rfloor-1}
\operatorname{osc}_{x\in M}
\left(
\partial_t v\left(
t-\bigl(\lfloor t/\tau\rfloor-1\bigr)\tau,\cdot
\right)
\right)\\
&\quad\leq
\theta^{\lfloor t/\tau\rfloor-1}
\operatorname{osc}_{x\in M}
\bigl(\partial_t v(0,\cdot)\bigr)\\
&\quad\leq
\theta^{\,t/\tau-2}
\operatorname{osc}_{x\in M}
\bigl(\partial_t v(0,\cdot)\bigr).
\end{aligned}
\end{equation*}
Here
\begin{equation*}
\tau
\leq
t-\bigl(\lfloor t/\tau\rfloor-1\bigr)\tau
<2\tau,
\end{equation*}
so every application of the contraction inequality starts at a positive time. The second inequality follows from Lemma \ref{thm:estimate of t derivative}, and the last follows from $\lfloor t/\tau\rfloor-1\geq t/\tau-2$ and $0<\theta<1$.

Since
\begin{equation*}
\int_M\partial_t v(t,\cdot)\,\omega^n=0,
\end{equation*}
we have
\begin{equation*}
-\operatorname{osc}_{x\in M}(\partial_t v)
\leq
\min_{x\in M}\partial_t v
\leq0
\leq
\max_{x\in M}\partial_t v
\leq
\operatorname{osc}_{x\in M}(\partial_t v).
\end{equation*}
Consequently,
\begin{equation*}
\sup_{x\in M}|\partial_t v(t,\cdot)|
\leq
\operatorname{osc}_{x\in M}
\bigl(\partial_t v(t,\cdot)\bigr).
\end{equation*}
Finally, the flow equation at $t=0$ and the spatial constancy of $c(0)$ imply
\begin{equation*}
\operatorname{osc}_{x\in M}
\bigl(\partial_t v(0,\cdot)\bigr)
=
\operatorname{osc}_{x\in M}
\Bigl(
\Delta v_0+\psi(|\nabla v_0|^2)
+\langle B,dv_0\rangle-f
\Bigr).
\end{equation*}
Therefore,
\begin{equation*}
\sup_{x\in M}|\partial_t v(t,\cdot)|
\leq
\operatorname{osc}_{x\in M}
\Bigl(
\Delta v_0+\psi(|\nabla v_0|^2)
+\langle B,dv_0\rangle-f
\Bigr)
\theta^{\,t/\tau-2},
\end{equation*}
which proves (\ref{eq:time derivative sup control}).
\end{proof}

We are now in a position to prove the convergence of the flow. The main result of this section is the following.

\begin{theorem}
As $t\to+\infty$, the solution $v(t,\cdot)$ converges uniformly on $M$ to a smooth function $v_\infty$ satisfying
\[
\Delta v_\infty+\psi(|\nabla v_\infty|^2)+\langle B,dv_\infty\rangle=f+c_\infty\qquad\text{on }M,
\]
where
    \begin{equation*}
        c_{\infty}
        =\frac{\int_{M}\big(\Delta v_{\infty}+\psi(\vert\nabla v_{\infty}\vert^2)
        +\langle B,dv_{\infty}\rangle\big)\,\omega^n}
        {\int_{M}\omega^n}.
    \end{equation*}
\end{theorem}

\begin{proof}
By Theorem \ref{thm:time derivative control}, there exist constants $C,\lambda>0$ and $t_0>0$ such that
\begin{equation*}
    \|\partial_t v(t,\cdot)\|_{L^\infty(M)}\le Ce^{-\lambda t}\qquad\text{for all }t\ge t_0.
\end{equation*}
Consequently, for $s>t\ge t_0$,
\begin{equation*}
   \begin{aligned}
        \|v(s,\cdot)-v(t,\cdot)\|_{L^\infty(M)}
        &\le\int_t^s\|\partial_r v(r,\cdot)\|_{L^\infty(M)}\,dr \\
         &\le \frac{C}{\lambda}e^{-\lambda t}.
    \end{aligned}
\end{equation*}
Thus there exists $v_\infty\in C(M)$ such that $v(t,\cdot)\to v_\infty$ uniformly on $M$ as $t\to+\infty$. It remains to show that $v_\infty$ is smooth and satisfies the stationary equation.

We first establish a uniform elliptic estimate for the time slices. At a maximum and a minimum point of $v(t,\cdot)$, respectively, the flow equation gives
\begin{equation*}
    |c(t)-\psi(0)|\le\|f\|_{L^\infty(M)}+\|\partial_t v(t,\cdot)\|_{L^\infty(M)}.
\end{equation*}
In particular, $c(t)$ is uniformly bounded for $t\ge t_0$.
Rewriting the flow equation as
\begin{equation*}
    \Delta v(t,\cdot)=\partial_t v(t,\cdot)-\psi(|\nabla v(t,\cdot)|^2)-\langle B,dv(t,\cdot)\rangle+c(t)+f,
\end{equation*}
and using Lemma \ref{thm:C_1 estimate} and Theorem \ref{thm:time derivative control}, we obtain
\begin{equation*}
    \sup_{t\ge t_0}\|\Delta v(t,\cdot)\|_{L^\infty(M)}\le C.
\end{equation*}
Since $\Delta$ is a fixed uniformly elliptic operator with smooth coefficients, the global elliptic $W^{2,p}$ estimate and Corollary \ref{thm:C_0 estimate} imply that, for every $1<p<\infty$,
\begin{equation*}
    \begin{aligned}
        \|v(t,\cdot)\|_{W^{2,p}(M)}&\le C_p\Bigl(\|\Delta v(t,\cdot)\|_{L^p(M)}+\|v(t,\cdot)\|_{L^p(M)}\Bigr) \\
         &\le C_p,\qquad t\ge t_0,
    \end{aligned}
\end{equation*}
where $C_p$ is independent of $t$.

Fix $p>2n$ and $0<\alpha<1-2n/p$, where $2n$ is the real dimension of $M$. Let $t_j\to+\infty$ and set $v_j=v(t_j,\cdot)$. By weak compactness in $W^{2,p}(M)$, the compact embedding
\begin{equation*}
    W^{2,p}(M)\Subset C^{1,\alpha}(M),
\end{equation*}
and the boundedness of $c(t_j)$, we may pass to a subsequence
such that
\begin{equation*}
    \begin{gathered}
        v_j\rightharpoonup v_\infty\quad\text{in }W^{2,p}(M),\\
        v_j\longrightarrow v_\infty\quad\text{in }C^{1,\alpha}(M),\\
        c(t_j)\longrightarrow c_\infty
    \end{gathered}
\end{equation*}
for some constant $c_\infty$. The function limit is $v_\infty$ because the entire flow already converges uniformly to $v_\infty$.

In particular, the strong convergence of the gradients yields
\begin{equation*}
    \psi(|\nabla v_j|^2)\longrightarrow\psi(|\nabla v_\infty|^2),\qquad\langle B,dv_j\rangle\longrightarrow\langle B,dv_\infty\rangle
\end{equation*}
uniformly on $M$. Since
\begin{equation*}
    \|\partial_t v(t_j,\cdot)\|_{L^\infty(M)}\longrightarrow0,
\end{equation*}
we can pass to the limit in the flow equation in the sense
of distributions to obtain
\begin{equation*}
    \Delta v_\infty+\psi(|\nabla v_\infty|^2)+\langle B,dv_\infty\rangle=f+c_\infty.
\end{equation*}
The right-hand side of
    \begin{equation*}
        \Delta v_\infty=f+c_\infty-\psi(|\nabla v_\infty|^2)-\langle B,dv_\infty\rangle
    \end{equation*}
belongs to $C^{0,\alpha}(M)$. Elliptic Schauder regularity therefore gives $v_\infty\in C^{2,\alpha}(M)$. Since $\psi$, $B$, $f$, and the coefficients of $\Delta$ are smooth, elliptic bootstrapping yields $v_\infty\in C^\infty(M)$.

Finally, integrating the stationary equation over $M$ and using $\int_M f\,\omega^n=0$, we obtain
    \begin{equation*}
        c_{\infty}
        =\frac{\int_{M}\big(\Delta v_{\infty}+\psi(\vert\nabla v_{\infty}\vert^2)
        +\langle B,dv_{\infty}\rangle\big)\,\omega^n}
        {\int_{M}\omega^n}.
    \end{equation*}
This completes the proof.
\end{proof}

\section{Uniqueness and generalized Gauduchon metrics}\label{sec:metrics}

It remains to prove the uniqueness in Theorem \ref{thm:semilinear eq}. 
Suppose that there exist $v$ and $\widetilde{v}$ satisfying 
\eqref{eq:semilinear eq} with
\begin{equation*}
    \Delta v+\psi(\vert\nabla v\vert^{2})+\langle B,dv\rangle =f+c,
\end{equation*}
and
\begin{equation*}
    \Delta\widetilde{v}+\psi(\vert\nabla\widetilde{v}\vert^{2})
    +\langle B,d\widetilde{v}\rangle =f+\widetilde{c}.
\end{equation*}
respectively. 
Here
\begin{equation*}
    c=\frac{\int_M(\Delta v+\psi(\vert\nabla v\vert^{2})
    +\langle B,dv\rangle)\omega^{n}}{\int_M\omega^{n}}
\end{equation*}
and
\begin{equation*}
    \widetilde{c}=\frac{\int_M(\Delta\widetilde{v}
    +\psi(\vert\nabla\widetilde{v}\vert^{2})
    +\langle B,d\widetilde{v}\rangle)\omega^{n}}
    {\int_M\omega^{n}}.
\end{equation*}
Thus, 
\begin{equation}\label{eq:uniquness of v}
    \Delta(v-\widetilde{v})
    +\Big(\psi(\vert\nabla v\vert^{2})
    -\psi(\vert\nabla\widetilde{v}\vert^{2})\Big)
    +\big\langle B,\,d(v-\widetilde{v})\big\rangle
    =c-\widetilde{c}.
\end{equation}
Applying maximum principle to $v-\widetilde{v}$ in (\ref{eq:uniquness of v}), we have $c=\widetilde{c}$.
Furthermore, using the strong maximum principle, we get that $v-\widetilde{v}$ is equal to a constant. The proof of Theorem \ref{thm:semilinear eq} is complete.

Set
\begin{equation}\label{eq:definition of B1}
    B_2=\frac{2k}{n-1}B_1
\end{equation}
and 
\begin{equation}\label{eq:definition of varphi}
    \varphi=\frac{n(\sqrt{-1}/2)\partial\bar\partial(\omega^k)\wedge\omega^{n-k-1}}{\omega^n}.
\end{equation}
We can rewrite (\ref{eq:main}) as
    \begin{equation*}
        \Delta v+\vert\nabla v\vert^2+\langle B_2,dv\rangle+\varphi=n\gamma_k(\omega).
    \end{equation*}
Then Theorem \ref{thm:main} follows immediately from Theorem \ref{thm:semilinear eq} with
\begin{equation*}
    \psi(t)=t\qquad \text{and} \qquad f=\frac{\int_M\varphi\omega^n}{\int_M\omega^n}-\varphi.
\end{equation*}
The corresponding geometric constant is
\begin{equation*}
    \gamma_k(\omega)
    =
    \frac{1}{n}
    \left(
        c+\frac{\int_M\varphi\,\omega^n}{\int_M\omega^n}
    \right).
\end{equation*}

\section*{Acknowledgements}

The author sincerely thanks Prof. Xi Zhang for the valuable guidance and insightful suggestions throughout this work.

\section*{Statements and Declarations}

The datasets used or analysed during the current study are available from the corresponding author on reasonable request.

\end{document}